\documentclass{article}
\usepackage{graphicx,amsmath,amssymb,amsthm} 
\usepackage{xcolor,float}

\title{Warping labeling for twisted knots and twisted virtual braids}
\author{Komal Negi, Ayaka Shimizu, Madeti Prabhakar }
\date{}
\newtheorem{theorem}{Theorem}[section]

\newtheorem{prop}[theorem]{Proposition}
\newtheorem{corollary}[theorem]{Corollary}

\theoremstyle{definition}
\newtheorem{definition}[theorem]{Definition}
\newtheorem{remark}[theorem]{Remark}
\newtheorem{Example}[theorem]{Example}
\begin{document}

\maketitle
\begin{abstract}
In this paper, we introduce the concept of the warping degree for twisted knots, construct an invariant for them, and utilize it to establish a labeling scheme for these knots, known as ``warping labeling". We have identified that a warping labeling can be extended to twisted virtual braids, enabling the creation of a function that remains invariant under all R-moves except the R2 move. By limiting the labeling set to $\mathbb{Z}_2$, we can develop invariants for twisted virtual braids.\\


\end{abstract}

\section{Introduction}

In this paper, twisted knots and twisted virtual braid diagrams are assumed to be oriented, and the diagrams are represented on the surface $S^2$.
The concept of twisted virtual braids was introduced in paper~\cite{KPS}. 
This work established significant results in the realm of braids, specifically proving Markov's and Alexander's theorems for twisted virtual braids and twisted links.
Currently, there are no known invariants for twisted virtual braids. 
In this work, we aim to define the first invariant for twisted virtual braids using the concept of warping degree. Many articles have explored the concept of warping degree in various contexts, including classical knots, virtual knots, and welded knots (\cite{IKL, LLW, KAY, A}).
The warping degree has been employed to develop invariants and establish lower bounds for unknotting operations in these knot categories.
In our research, we extend the notion of warping degree to twisted knots and further construct labeling for twisted knots and twisted virtual braids. 
This extension allows us to define new invariants that can be used to characterize and distinguish twisted virtual braids.

The structure of the article is outlined as follows.
Section 2, provides essential background information on twisted knots, twisted virtual braids, and the concept of warping degree for knots, setting the foundation for our research.
In Section 3, we introduce the concept of warping degree specifically tailored for twisted knot diagrams, expanding on the traditional understanding to accommodate the complexities of twisted structures. Also, we construct an invariant for twisted knots using warping degree.
In Section 4, we define a labeling method for both twisted knot diagrams and twisted virtual braid diagrams. The labeling for twisted virtual braids helps us identify the presence of R2 moves between two equivalent twisted virtual braids.
Section 5, introduces a novel labeling approach for twisted virtual braids known as $\mathbb{Z}_2$-labeling, and establishes its invariance property for twisted virtual braids, providing a new tool for studying and classifying these structures.
The conclusion of the article highlights open problems and future directions related to the topics discussed, suggesting avenues for further research and exploration in the field of twisted knots and twisted virtual braids.
\section{Preliminaries}
Twisted knots are stable equivalence classes of oriented knots in orientable three-manifolds that are orientation I-bundles over closed but not necessarily orientable surfaces. 
 
 There is a bijection between ambient isotopy equivalence classes of knots in stable oriented thickenings and Reidemeister equivalence classes of twisted knot diagrams. We have taken the preliminaries from the papers~\cite{MO, KPS, A}.
\subsection{Twisted knot diagrams}
\begin{figure}[h]
	\centering
\includegraphics[width=9cm,height=5cm]{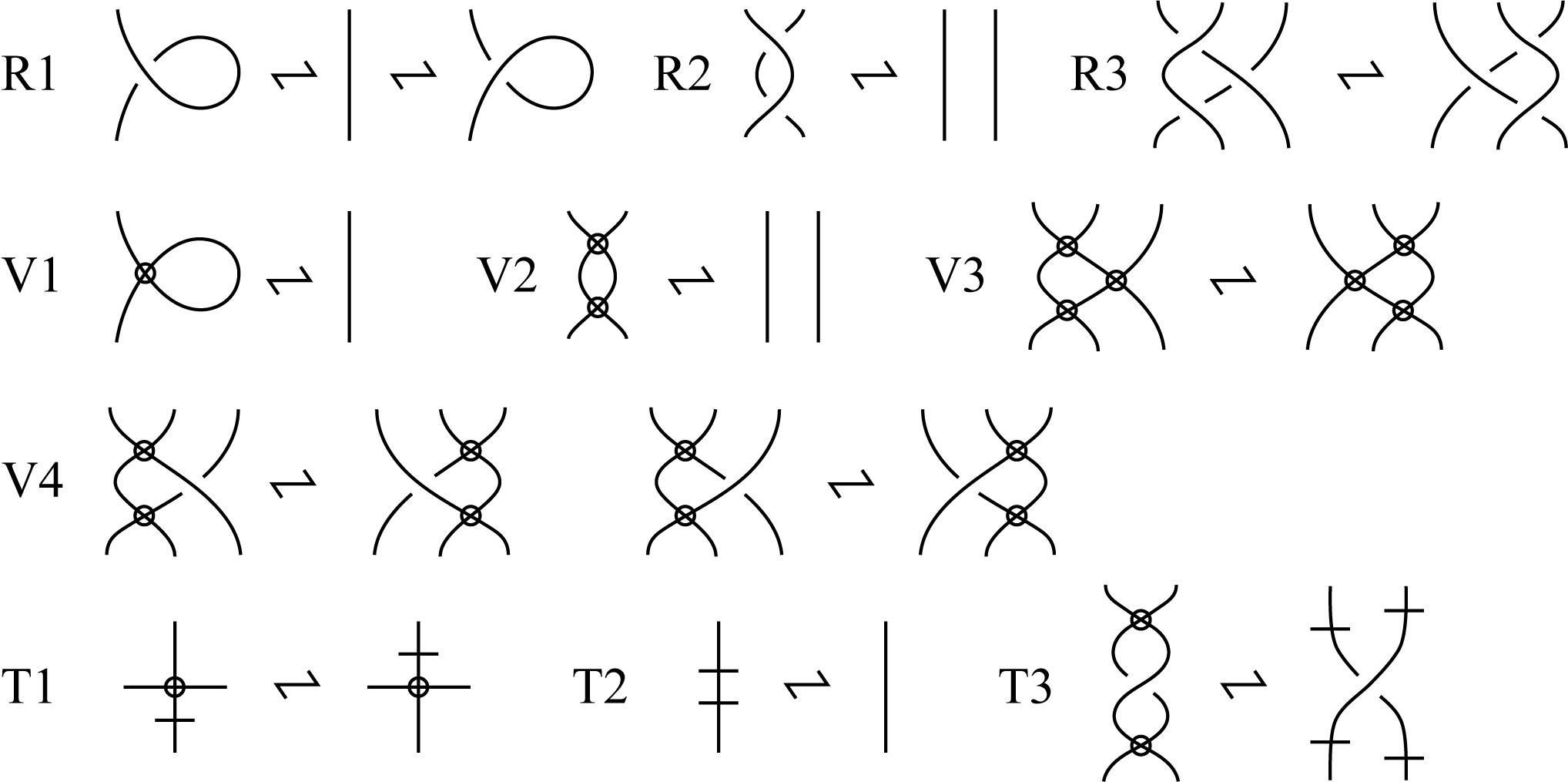}	
  \caption{Extended Reidemeister moves for twisted knot diagrams}
 \label{rt}
   \end{figure}

\textit{Twisted knot diagrams}~\cite{MO} are defined as marked generic planar curves, where the markings identify the usual classical crossings, virtual crossings, and bars on edges. Reidemeister moves for twisted link theory is shown in Figure~\ref{rt}. Classical Reidemeister moves (R1, R2, R3), virtual Reidemeister moves (V1, V2 V3,V4) and twisted Reidemeister moves (T1, T2, T3), are together called as extended Reidemeister moves.

It is known that two diagrams represent ambient isotopic twisted diagrams if and only if they are transformed into
each other by a ﬁnite sequence of extended Reidemeister moves, as shown in Figure~\ref{rt}. 

\begin{definition}
    \textit{Trivial twisted knot} is either a trivial knot or a trivial knot with a bar.
\end{definition}

\subsection{Twisted virtual braids}
Let $n$ be a positive integer.  

\begin{figure}[h]
  \centering
    \includegraphics[width=2cm,height=3cm]{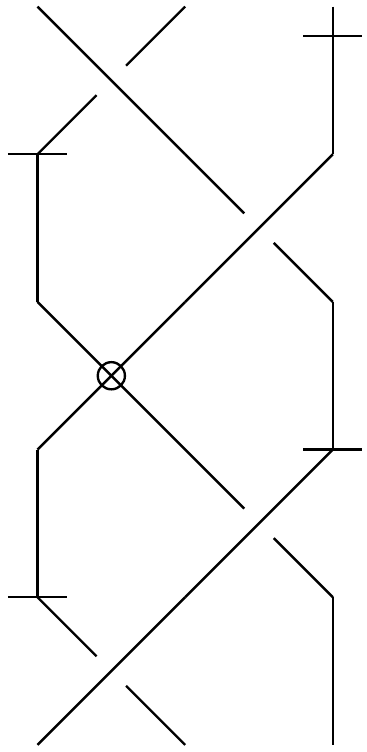}
        \caption{A twisted virtual braid on 3 strands}
        \label{exa}
        \end{figure}

\begin{definition}
A {\it twisted virtual braid diagram} on $n$ strands (or of degree $n$) 
is a union of $n$ smooth or polygonal curves, which are called {\it strands}, in $\mathbb{R}^2$ connecting points $(i,1)$ with points $(q_i,0)$ $(i=1, \dots, n)$, where $(q_1, \ldots, q_n)$ is a permutation of the numbers $(1, \ldots, n)$, such that these curves are monotonic with respect to the second coordinate and intersections of the curves are transverse double points equipped with information as a positive/negative/virtual crossing and 
strings may have {\it bars} by which we mean short arcs intersecting the strings transversely.
See Figure~\ref{exa}, where the five crossings are 
negative, positive, virtual, positive and positive from the top. 
\end{definition}

\begin{figure}[h]
  \centering
    \includegraphics[width=8cm]{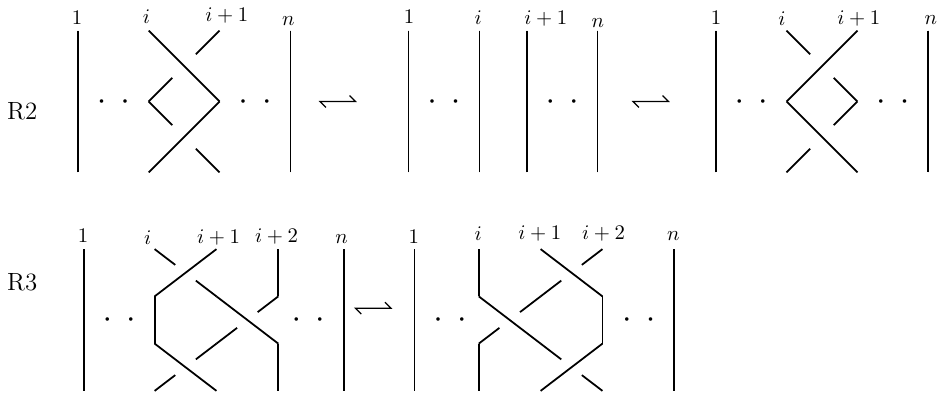}
        \caption{Classical braid moves}
        \label{bmoves}
        \end{figure}

\begin{figure}[h]
  \centering
    \includegraphics[width=10cm]{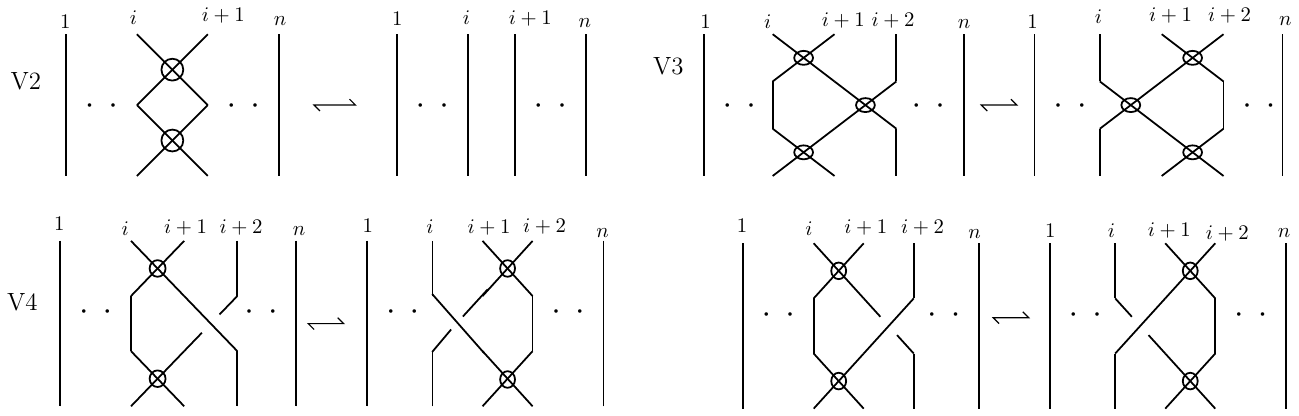}
        \caption{Virtual braid moves}
        \label{vbmoves}
        \end{figure}

\begin{figure}[h]
  \centering
    \includegraphics[width=9cm]{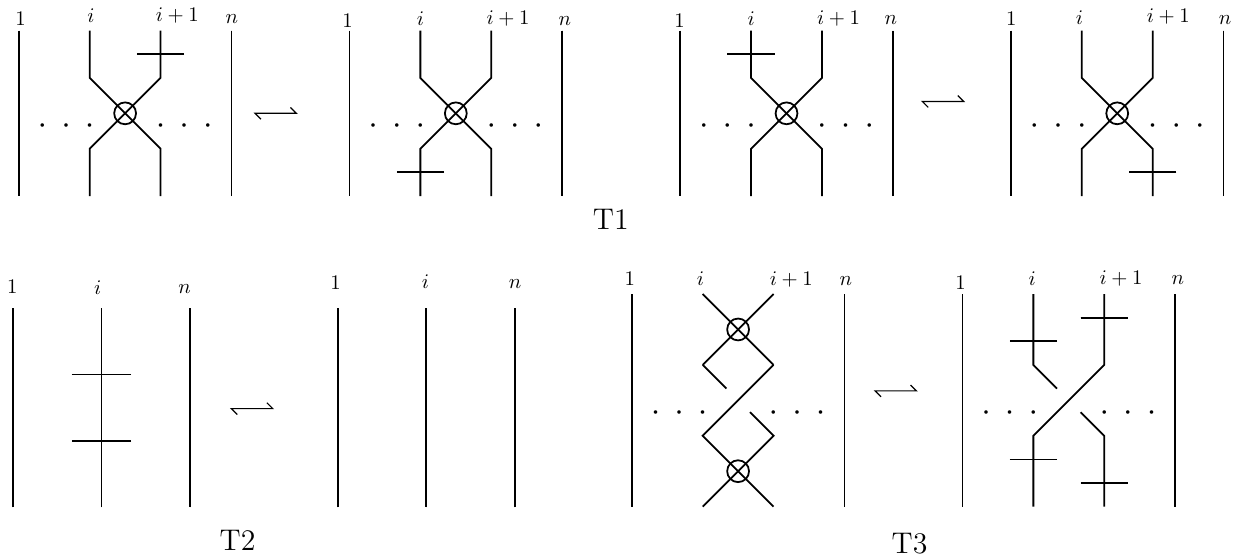}
        \caption{Twisted braid moves}
        \label{moves}
        \end{figure}

\begin{definition}
Two twisted virtual braid diagrams $\beta$ and $\beta'$ of degree $n$ are {\it equivalent} if they are related by classical, virtual, and twisted braid moves shown in Figure~\ref{bmoves}, \ref{vbmoves}, and \ref{moves}, respectively. Collectively, these are referred to as R-moves.
\end{definition}

\begin{definition}
     A {\it twisted virtual braid} is an equivalence class of twisted virtual braid diagrams. 
\end{definition}

\subsection{Warping degree of knot diagrams}
Consider $D$ as an oriented knot diagram or virtual knot diagram. Let $a$ be a point on  $D$  that is not a crossing point; we refer to $a$ as a base point of $D$. The pair $(D, a)$ is denoted as $D_a$. A crossing point of $D_a$ is termed a \textit{warping crossing} point if, when traversing along the oriented diagram $D$ starting from $a$, we encounter the crossing point first as an under-crossing.
\begin{definition}
The \textit{warping degree} of $ D_a$, represented as $d(D_a)$, is defined as the count of warping crossing points in $D_a$. The warping degree of  $D$, denoted as $d(D)$, is the smallest warping degree observed across all base points of $D$.
\end{definition}

\noindent By definition, warping degree represents how far a diagram is from a descending diagram with given orientation. 

\begin{Example}
    Consider the knot diagram shown in Figure~\ref{exk}. The warping degree of the diagram is $1$.
\begin{figure}[h]
	\centering
\includegraphics[width=2cm]{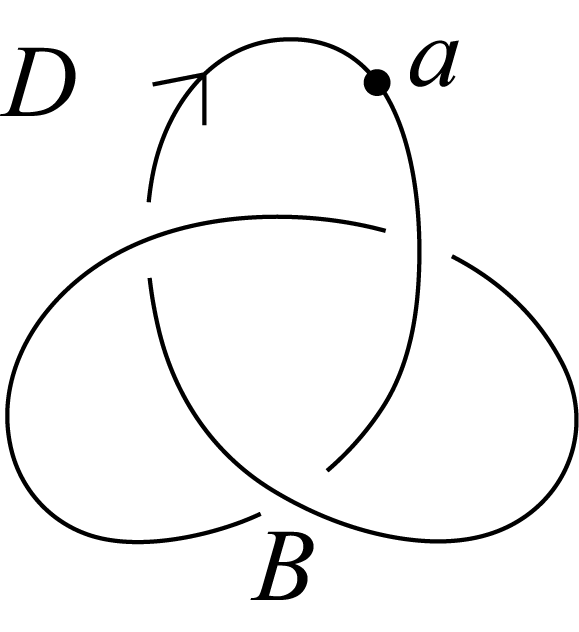}	
  \caption{The crossing $B$ is a warping crossing point of $D_a$}
 \label{exk}
   \end{figure}
\end{Example}
\begin{prop}[\cite{A}]
    Let $D$ be an oriented knot diagram with crossing number $\# C(D)$, and let $-D$ denote $D$ with orientation reversed. Then $d(D_a)+d(-D_a)= \# C(D)$ holds for any base point $a$ on $D$. 
\label{dada}
\end{prop}
\begin{theorem}[\cite{A}]
Let D be an oriented knot diagram which has at least one crossing point. Then we have the following inequality:
$$d(D) + d(-D) + 1 \leq \# C(D).$$
Further, the equality holds if and only if D is an alternating diagram.
\label{d-d}
\end{theorem}
    
Motivated by the above work we define warping degree for twisted knot diagrams.

\section{Warping degree of twisted knot diagrams}
\begin{definition}
Consider $D_a$ as an oriented twisted knot diagram with  $a$ as a base point of $D$.  A crossing $C$ of $D_a$ is called a \textit{warping crossing} if it satisfies either of the following condition:
\begin{itemize}
    \item  When traversing along the oriented diagram $D$ starting from $a$, we encounter the crossing point first as an under-crossing and there is an even number of bars from the base point $a$ to the crossing $C$. 
    \item When traversing along the oriented diagram $D$ starting from $a$, we encounter the crossing point first as an over-crossing and there is an odd number of bars from the base point $a$ to the crossing $C$.
\end{itemize}
The \textit{warping degree} of $ D_a$, is defined as the number of warping crossings in $D_a$.
\end{definition}
\begin{definition}
For an oriented twisted knot diagram $D$, the \textit{warping degree} of $D$ is $d(D)=\min\{D_a| \ a \text{ is a base point on } D\}$.
\end{definition}
    \begin{Example}
    Consider a twisted knot diagram shown in Figure~\ref{ext}. Observe that, $d(D_a)=1$, $d(D_b)=0$ and $d(D_c)=2$. Also, the warping degree of the given twisted knot diagram is $0$. We note that $d(D_a)+d(-D_a)=1 \neq 2=\#C(D)$, that is Proposition \ref{dada} does not hold for the twisted knot diagram. 
    \begin{figure}[h]
	\centering
\includegraphics[width=2cm]{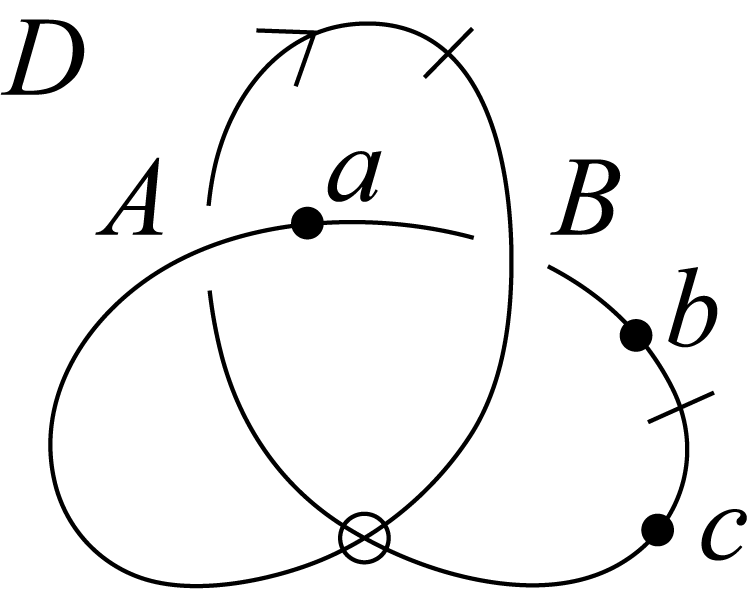}	
  \caption{A twisted knot diagram $D$ with base points $a$, $b$ and $c$}
 \label{ext}
   \end{figure}
\end{Example}

\begin{prop}\label{RWT}
           Let $a$ and $b$ be the base points on the diagram $D$ located on opposite sides of a classical crossing or a bar as shown in Figure \ref{pbp}. Let $\#C(D)$ be the number of classical crossings of the diagram $D$. Then the following statements are true:
           \begin{itemize}
               \item[1.] If there exists a bar between $a$ and $b$ then $d(D_a)+d(D_b)=\#C(D)$.
               \item[2.] If there is a crossing $B$ between $a$ and $b$ then the following happens.
               \begin{itemize}
                   \item[(i)] If there are odd bars in the path from $B$ to $B$ which includes $b$, then $d(D_a)=d(D_b)$.
                   \item[(ii)] If there are even bars from $B$ to $B$ which includes $b$ and $B$ is an over-crossing (resp. under-crossing), then $d(D_b)=d(D_a)+1$ (resp. $d(D_b)=d(D_a)-1$).
               \end{itemize}
           \end{itemize}
\end{prop}
\begin{proof}Let $a$ and $b$ be the base points on the diagram $D$. Let $\# C(D)$ be the number of crossings of the diagram $D$.
\begin{itemize}
    \item[1.] If there is a bar between $a$ and $b$ then all the non-warping crossing of $D_a$ are the warping crossings of $D_b$. Observe that the sum of the number of warping crossings of $D_a$ and the number of non-warping crossings of $D_a$ gives $\# C(D)$.

    \item[2.] If there is a crossing $B$ between $a$ and $b$ and odd bars presents from $B$ to $B$ moving along the orientation as shown in Figure~\ref{pbp}. Suppose $C$ is a crossing such that $C\neq B$, then nature of $C$ will be same for $D_a$ and $D_b$. If $C=B$ and $C$ is a warping crossing for $D_a$ then $C$ will be warping crossing for $D_b$ due to the odd number of bars from $b$ to $C$, and vice-versa. Therefore, $d(D_a)=d(D_b)$.  
  \begin{figure}[h]
  \centering
    \includegraphics[width=6.5cm]{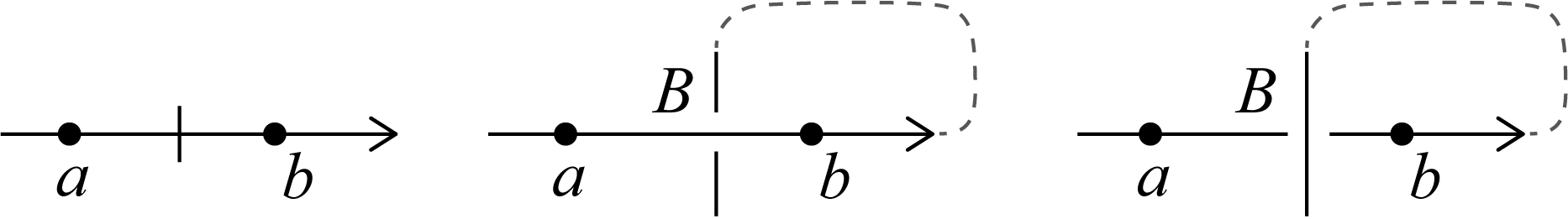}
        \caption{Different positions of base points on a diagram}
        \label{pbp}
        \end{figure}
    If there are even bars from $B$ to $B$, then we can proceed similarly as done for classical knot diagrams in~\cite{A}.
\end{itemize}
\end{proof}
\begin{prop}
    Let $D$ be an oriented twisted knot diagram. Let $d(D)$ denote the warping degree of $D$, and $\# C(D)$ the crossing number of $D$. The following relation holds, $$d(D)\leq \frac{\# C(D)}{2}.$$ 
\label{DC}
\end{prop}
\begin{proof}
    There exist two base points $a$ and $b$ such that there is a bar between $a$ and $b$, which satisfy the following,
    $$d(D_a)+d(D_b)=\#C(D).$$
    Since $d(D)=\min\{D_a| \ a \text{ is a base point on } D\}$, we have $d(D)\leq d(D_a)$ and $d(D)\leq d(D_b)$. Now the following follows,
     $d(D)+d(D)\leq d(D_a)+d(D_b)=\#C(D)$ or $2d(D)\leq \#C(D)$.
\end{proof}

\begin{remark}
Proposition \ref{DC} holds for any twisted knot diagram with bars for any orientation. For classical knot diagrams without bars, Proposition \ref{DC} does not hold. See Figure \ref{d2c3}. 
Indeed, $d(D) \leq ( \# C(D)-1)/2$ holds for classical knot diagrams if we take an appropriate orientation by Theorem \ref{d-d}. 
\begin{figure}[h]
	\centering
\includegraphics[width=2cm]{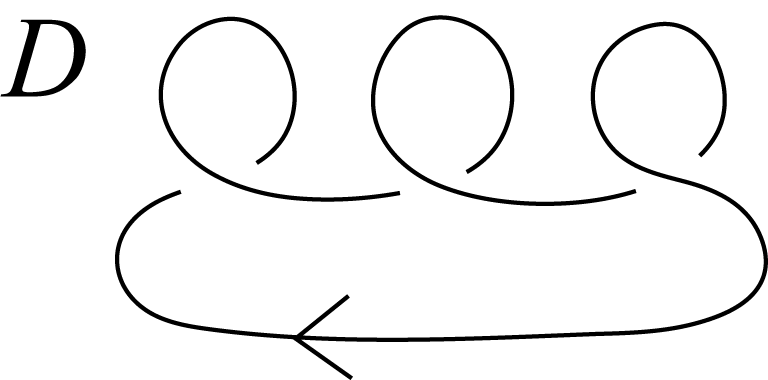}	
  \caption{An oriented diagram $D$ with $d(D)=2$ and $\# C(D)=3$}
 \label{d2c3}
\end{figure}
\end{remark}

\begin{definition}
    Let $D$ be a twisted knot diagram and $D^*$ is said to be mirror image of $D$ if $D^*$ is obtained from 
$D$ by changing the crossing at every crossing point.
\end{definition}
\begin{prop}
    Let $D_a$ be an oriented twisted knot diagram with base point $a$ and $D_a^*$ be a mirror image of $D_a$ with the same orientation. Then we have the following inequality
    $$d(D_a)+ d(D_a^*)=\#C(D).$$
\end{prop}
\begin{proof}
    The non-warping crossings of $D_a$ are the warping crossings of $D_a^*$.
\end{proof}
\begin{theorem}
     Let $D$ be an oriented twisted knot diagram and $D^*$ be a mirror image of $D$ with the same orientation. Then we have the following inequality
    $$d(D)+ d(D^*)\leq \#C(D).$$
\end{theorem}
\begin{definition}
   The \textit{warping degree of a twisted knot $K$} is $$d(K)=\min \{d(D),d(-D)| D \text{ represents } K\}.$$ It is an invariant for twisted knots.
\end{definition}

\section{Warping labeling}
In this section, we introduce a labeling process based on Proposition~\ref{RWT}, that allows us to label twisted knot diagrams.

\subsection{Warping labeling for twisted knots }
In this section, we refer to each part of an arc of a twisted knot diagram between classical crossings or bars as an edge. 
\begin{definition}
    The \textit{warping labeling} is a labeling giving each edge the value of warping degree with a base point on that edge. 
\end{definition}
\begin{figure}[h]
	\centering
\includegraphics[width=4.5cm]{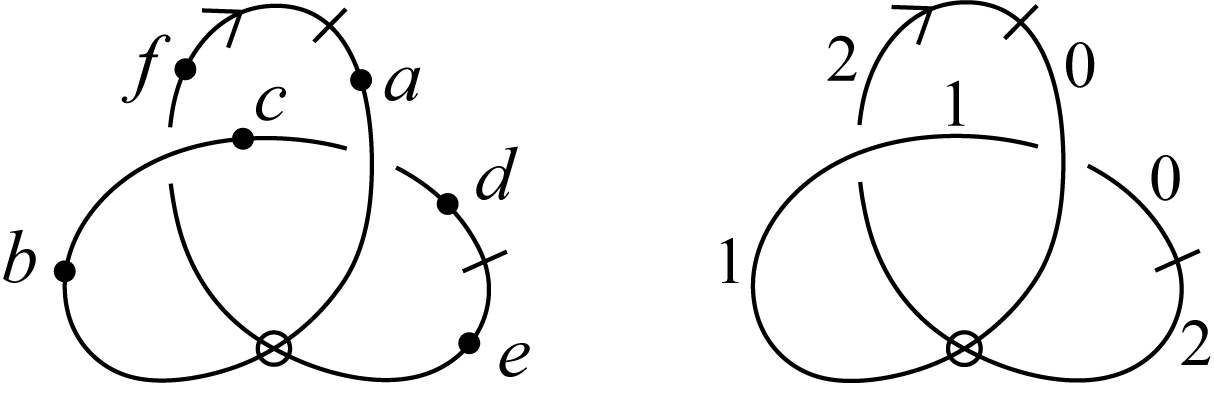}	
  \caption{The warping labeling}
 \label{wl}
\end{figure}

\begin{definition}
    An \textit{up-down labeling} is a labeling to each edge of an oriented twisted knot diagram which satisfies the rule of Proposition \ref{RWT}.
\end{definition}
\begin{figure}[h]
	\centering
\includegraphics[width=4.5cm]{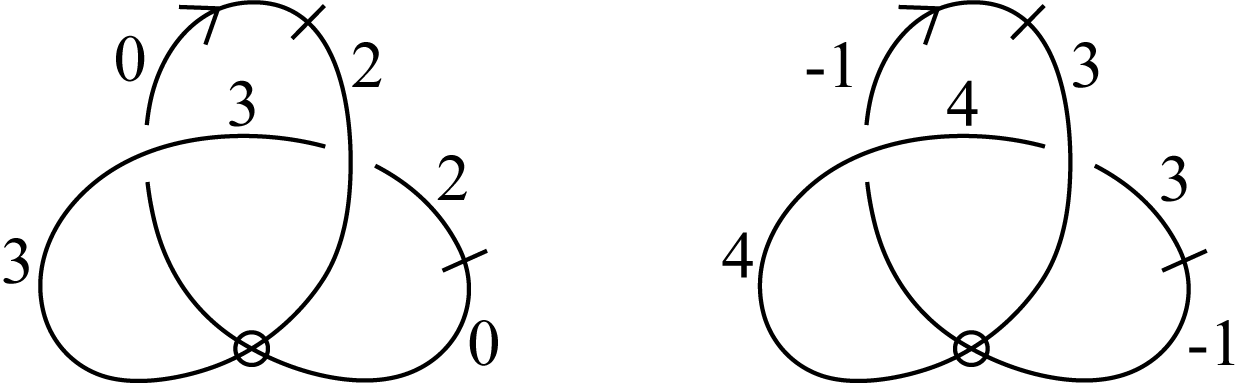}	
  \caption{Up-down labelings}
 \label{}
\end{figure}
\noindent The following proposition may be useful when giving an up-down labeling. 
\begin{prop} Let $D$ be an oriented twisted knot diagram with an up-down labeling. 
\begin{itemize}
    \item[(1)] When $D$ has an odd number of bars, the labeling around each crossing is as shown in Figure \ref{cross}, (1). 
    \item[(2)] When $D$ has an even number of bars, the labeling around each crossing is either of the two types shown in Figure \ref{cross}, (2).
\end{itemize}
\begin{figure}[h]
	\centering
\includegraphics[width=5cm]{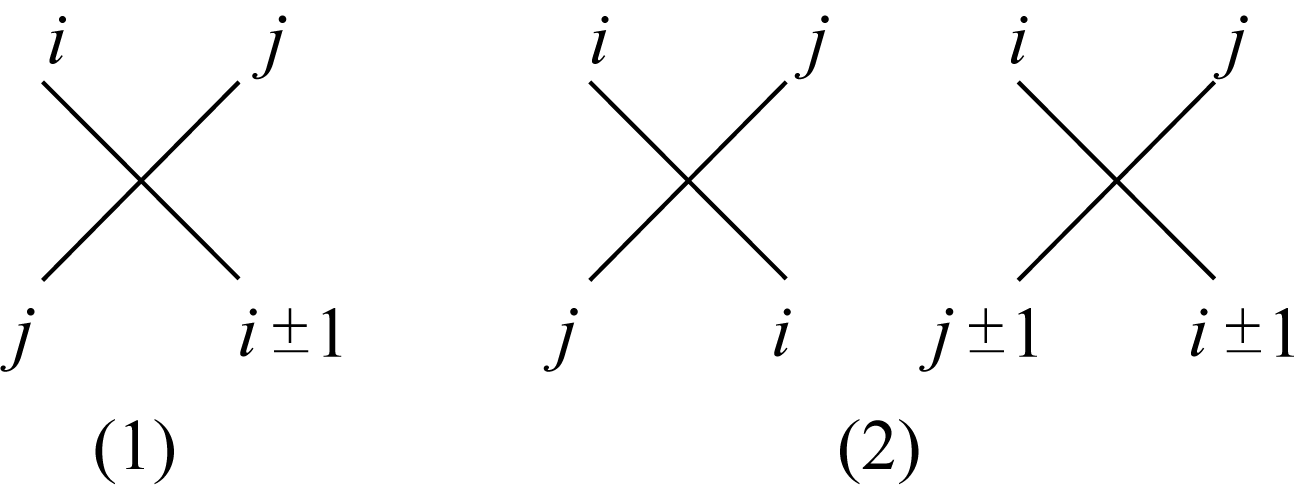}	
  \caption{(1) When $D$ has an odd number of bars, (2) When $D$ has an even number of bars }
 \label{cross}
\end{figure}
\end{prop}
\begin{proof}
    Let $D$ be an oriented twisted knot diagram. For each crossing $B$, $D$ is divided into two parts, the path `$C_1$' from $B$ to $B$ starting after the over-crossing of $B$, and the path `$C_2$' from $B$ to $B$ starting after the under-crossing of $B$, as shown in Figure \ref{cross2}. Since each bar belongs to just one of $C_1$ or $C_2$, the sum of the numbers of bars in $C_1$ and $C_2$ is equal to the number of bars of $D$. (1) When $D$ has an odd number of bars, one of $C_1$ or $C_2$ has an odd number of bars, and the other one has an even number of bars. (2) When $D$ has an even number of bars, then both of $C_1$ and $C_2$ have an even number of bars, or both have an odd number of bars. 
\begin{figure}[h]
	\centering
\includegraphics[width=2.5cm]{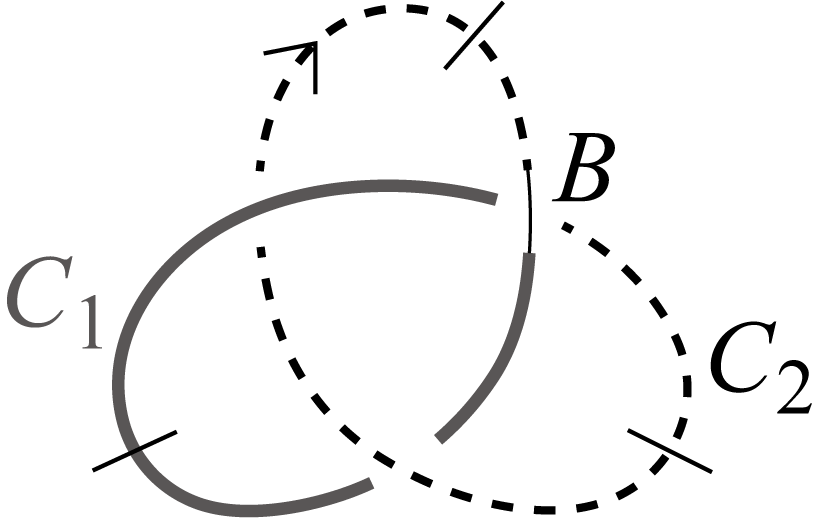}	
  \caption{A diagram divided at a crossing $B$}
 \label{cross2}
\end{figure}    
\end{proof}

\begin{figure}[h]
	\centering
\includegraphics[width=2.9cm]{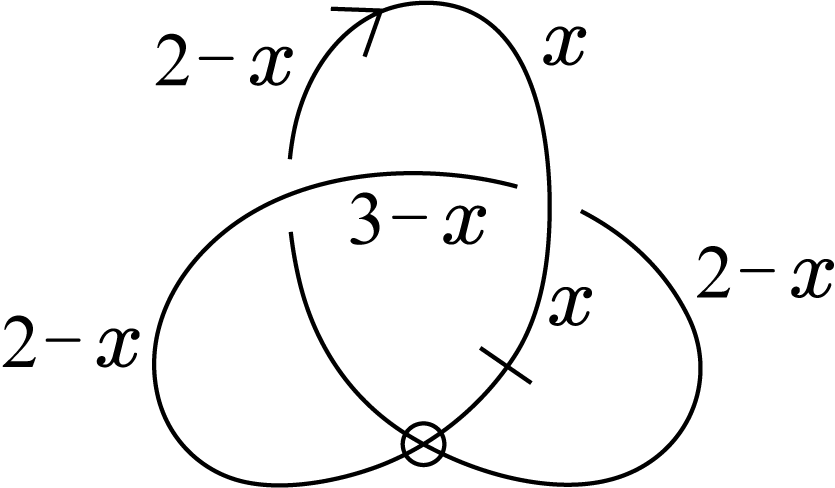}	
  \caption{By solving $2-x=x$, the value of $x$ must be $1$}
 \label{ul}
   \end{figure}

\noindent Number of the up-down labeling for twisted knots will be either unique or infinite. For instance, twisted knots with no bars, will have infinite number of such labeling, for details see~\cite{KAY}.
For some twisted knots this labeling is unique, see Figure~\ref{ul}. 

\begin{prop}
    An oriented twisted knot diagram $D$ has a unique up-down labeling if $D$ has an odd number of bars, and has infinitely many up-down labeling if $D$ has an even number of bars. 
\end{prop}
\begin{proof}
    Firstly, $D$ has at least one up-down labeling, which is the warping labeling. For an oriented knot diagram $D$, choose any edge $e$ to take a start point. Give a label $x$ to $e$. Then the label to the next edge is uniquely determined as a linear expression of $x$ according to the rule of Proposition \ref{RWT}. In the same way, all the edges receive linear expressions of $x$. Notice that the coefficient of $x$ is multiplied by $-1$ at every bar. After traversing $D$, the edge $e$ is again assigned a linear expression, $ax+b$. Notice that $a=-1$ if $D$ has an odd number of bars and $a=1$ if $D$ has an even number of bars. Comparing to the initial labeling $x$, we obtain a linear equation $ax+b=x$ in which all the solutions admit the up-down labeling. Recall that the equation $ax+b=x$ has at least one solution, the warping labeling. 
    If $D$ has an odd number of bars, the equation $-x+b=x$ has a unique solution. If $D$ has an even number of bars, then the equation $x+b=x$ should satisfy $b=0$, otherwise it has no solutions. Then the equation $x=x$ has infinitely many solutions.
\end{proof}

\begin{Example} 
Consider a family of twisted knot diagrams denoted as TD(n) such that $\#C(\text{TD(n)})=2n$ and each diagram contains 2 bars, as shown in Figure~\ref{fwl}. For the diagram TD(n), the warping labeling of each edge is equal to $n$. In this example, the inequality of Proposition~\ref{DC} is best possible, and also for any  number $n$, there exists a twisted knot diagram $D$ such that $ d(D)=n$.
    \begin{figure}[h]
	\centering
\includegraphics[width=6cm]{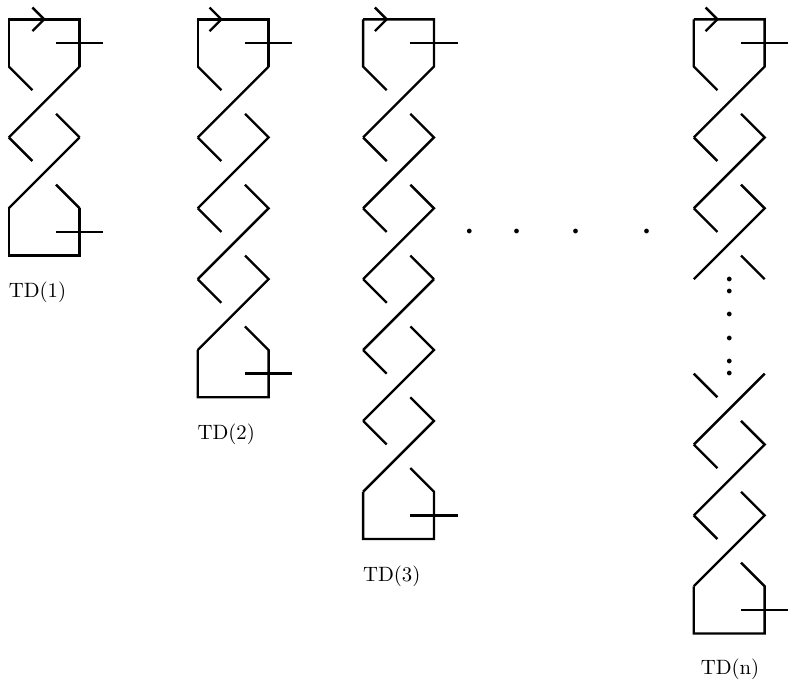}	
  \caption{A family of twisted knot diagrams TD(n) with warping degree $n$}
 \label{fwl}
   \end{figure}
\end{Example}

\subsection{Up-down labeling for twisted virtual braids}
The labeling method used for twisted knots can similarly be applied to label twisted virtual braids and later we will construct invariant for the set $TVB_n$. The method of labeling the twisted virtual braid is using Proposition~\ref{RWT}. Let $\beta$ be a twisted virtual braid diagram. We refer to each part of a strand divided by classical crossings or bars as an edge. 
\begin{definition}
    For a twisted virtual braid diagram $\beta$, a labeling to each edge is said to be an \textit{up-down labeling} to $\beta$ if the labeling satisfies the following conditions. Let $e_1$, $e_2$ be edges of $\beta$ located on opposite sides of a classical crossing or a bar as shown in Figure \ref{pbp2}. Let $l(e_1)$, $l(e_2)$ be the labels for $e_1, e_2$, respectively. Let $\#C(\beta)$ be the number of classical crossings of the braid diagram $\beta$. 
\begin{itemize}
\item[1.] If there exists a bar between $e_1$ and $e_2$ then $l(e_1)+l(e_2)=\#C(\beta)$.
\item[2.] If there is a crossing $B$ between $e_1$ and $e_2$ then the following happens.
\begin{itemize}
\item[(i)] If there are odd bars in the path from $B$ to $B$ on the closure of $\beta$ which includes $e_2$, then $l(e_1)=l(e_2)$.
\item[(ii)] If there are even bars from $B$ to $B$ and $B$ is an over-crossing (resp. under-crossing), then $l(e_2)=l(e_1)+1$ (resp. $l(e_2)=l(e_1)-1$).
\end{itemize}
\end{itemize}
\begin{figure}[h]
	\centering
\includegraphics[width=9cm]{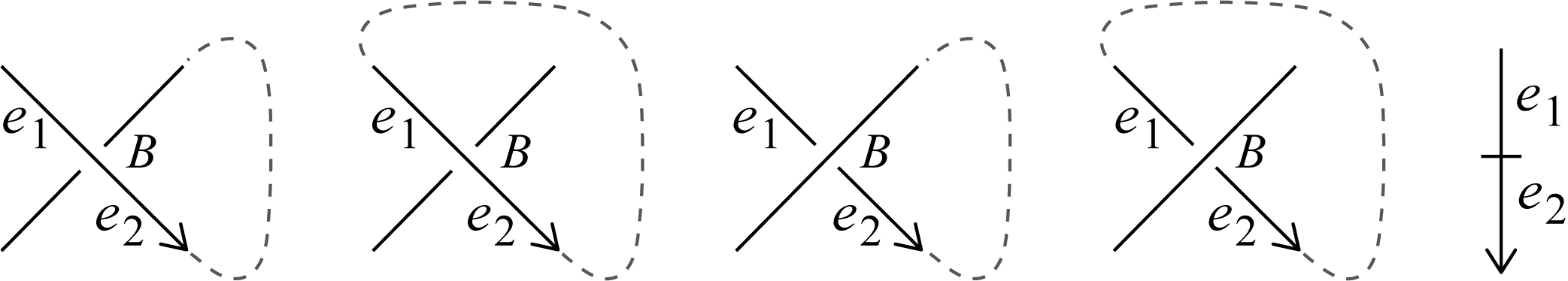}	
  \caption{Two edges $e_1$ and $e_2$}
 \label{pbp2}
\end{figure}
\end{definition}
\begin{Example}
    For the twisted virtual braid of three strands in Figure \ref{wdll}, give any number on the top for each strands. For example, we give $(1,2,3)$ here. Following the rule, give the up-down labeling. (1) On the first strand, the path from the crossing $c_1$ to $c_1$ in the closure includes two bars. Hence, the label is increased by one after the over-crossing $c_1$. (2) On the first strand, the path from $c_2$ to $c_2$ includes one bar. Hence the label is preserved at $c_2$ on the first strand. (3) The up-down labeling for the initial tuple $(1,2,3)$ can be given uniquely.
\begin{figure}[h]
	\centering
\includegraphics[width=9cm]{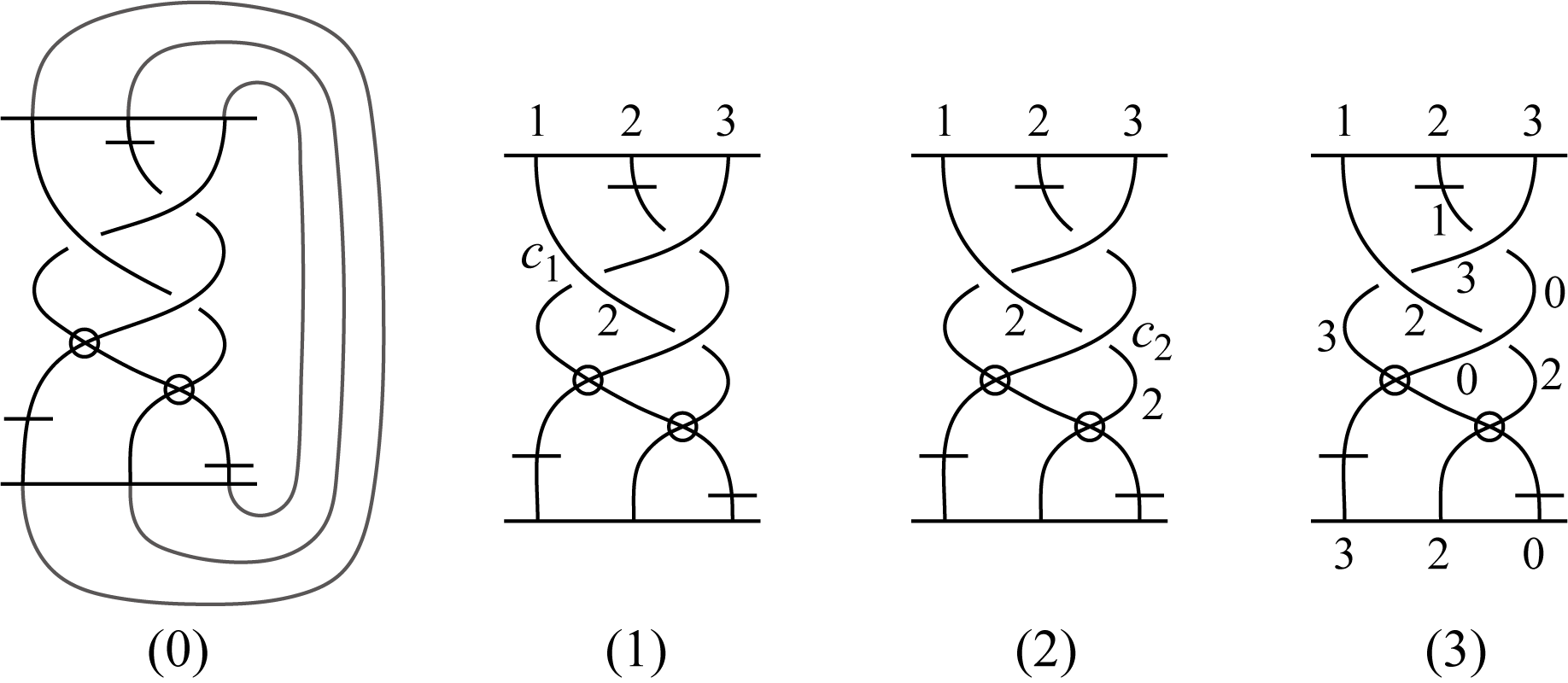}	
  \caption{An up-down labeling }
 \label{wdll}
\end{figure}
\end{Example}

\begin{definition}
    Let $\beta$ be a twisted virtual braid diagram, $f_\beta: \mathbb{Z}^n \to \mathbb{Z}^n$ such that $f_\beta(x_1,x_2,\ldots,x_n)=(y_1,y_2,\ldots,y_n)$, where $x_i$ is the initial labeling of the $i^{th}$ strand and $y_j$ is the labeling of $j^{th}$ strand on the bottom of $\beta$ for edges labeled an up-down labeling. The function $f_\beta$ is called \textit{up-down labeling function}.
\begin{figure}[h]
\centering
\includegraphics[width=2.5cm]{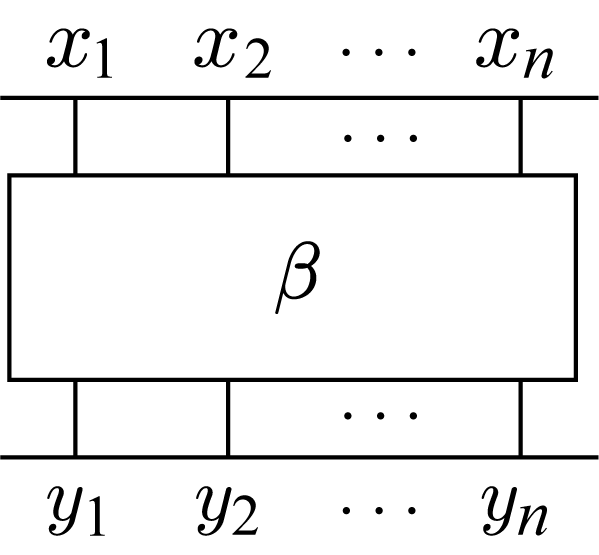}	
  \caption{The up-down labeling function $f_\beta$}
 \label{func}
\end{figure}
\end{definition}

\begin{Example}
The up-down labeling function $f_\beta$ with the braid $\beta$ in Figure~\ref{wdll2} has the outputs $f_\beta (1,2,3)=(3,2,0)$, $f_\beta (1,1,1)=(2,2,2)$. Moreover, $f_\beta(x_1,x_2,x_3)=(1+x_2,1+x_1,3-x_3)$.
\begin{figure}[h]
	\centering
\includegraphics[width=6cm]{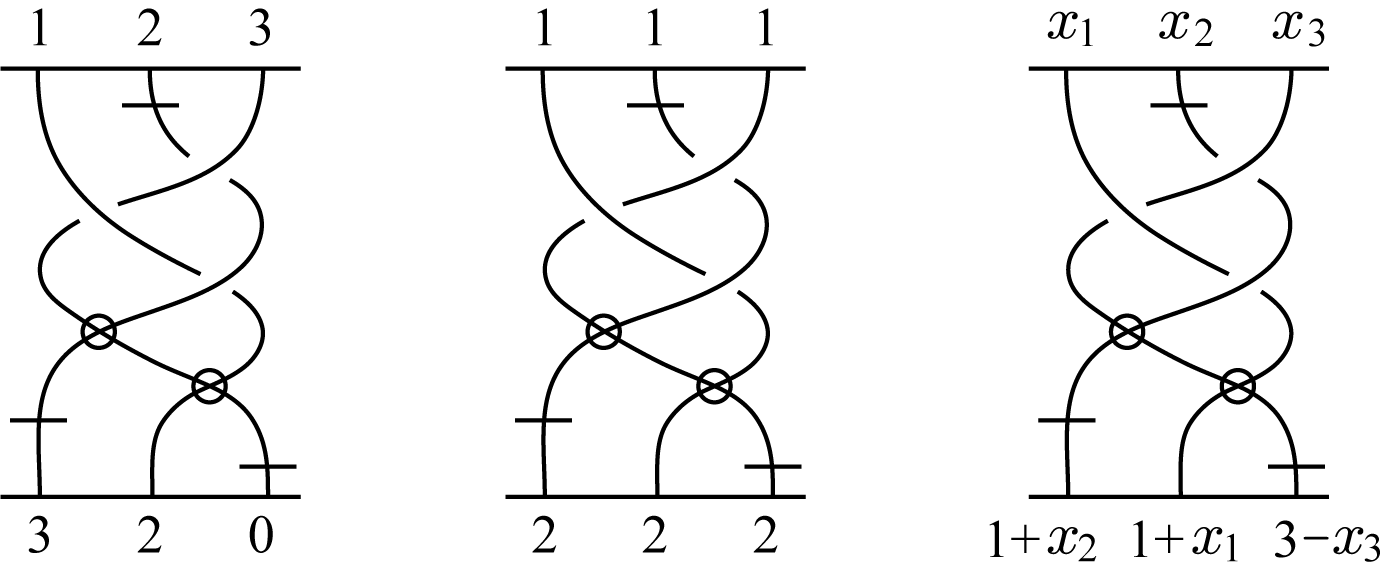}	
  \caption{$f_\beta(x_1,x_2,x_3)=(1+x_2,1+x_1,3-x_3)$}
 \label{wdll2}
\end{figure}
\end{Example}

\begin{remark}
    Let $\beta$ be a twisted virtual braid with $n$ strands and $p_\beta$ be the permutation of $\beta$. Let some linear function of variable $x$ is denoted by $L_x(x)$. Then the up-down labeling function $f_\beta(x_1,x_2,\ldots,x_n)=(y_1,y_2,\ldots,y_n)$ satisfy the following relation: $y_{p_{\beta}(i)}=L_{x_i}(x_i)$. 
\end{remark}
\begin{remark}
    If $\beta$ is a twisted virtual pure braid~\cite{VTKM}, then $f_\beta(x_1,x_2,\ldots,x_n)=(L_{x_1}(x_1),L_{x_2}(x_2),\ldots,L_{x_n}(x_n))$.
\end{remark}
\begin{theorem}\label{r2indicator1}
    If $\beta$ and $\beta'$ are equivalent twisted virtual braids upto all R-moves except R2, then the map $f_\beta(x_1,x_2,\ldots,x_n)=f_{\beta'}(x_1,x_2,\ldots,x_n)$, for any $(x_1,x_2,\ldots,x_n) \in \mathbb{Z}^n$.
\end{theorem}
\begin{proof}
It is sufficient to prove that the function $f_\beta$ is invariant under all R-moves except $R2$ move.

For the labeling under R3 move, see Figure \ref{r3}. The two parts of a braid are related by an R3 move. Give the same labels, $i, j, k$ on the top. 
There is a correspondence between the crossings $c_1$ and $d_1$ as the crossing between the first and second strands. After passing through the crossings $c_1$ and $d_1$ on the first strand, there are the same changes, say $\varepsilon_1^2$; if there are an odd number of bars from $c_1$ to $c_1$, then there are also an odd number of bars from $d_1$ to $d_1$ since the outside are same. In this case, $\varepsilon_1^2 =0$. If there are an even number of bars, $\varepsilon_1^2 =-1$. In the same way, there are the same changes at corresponding crossings, and as a result they have the same labeling on the bottom. 
\begin{figure}[H]
 \centering
 \includegraphics[width=7cm]{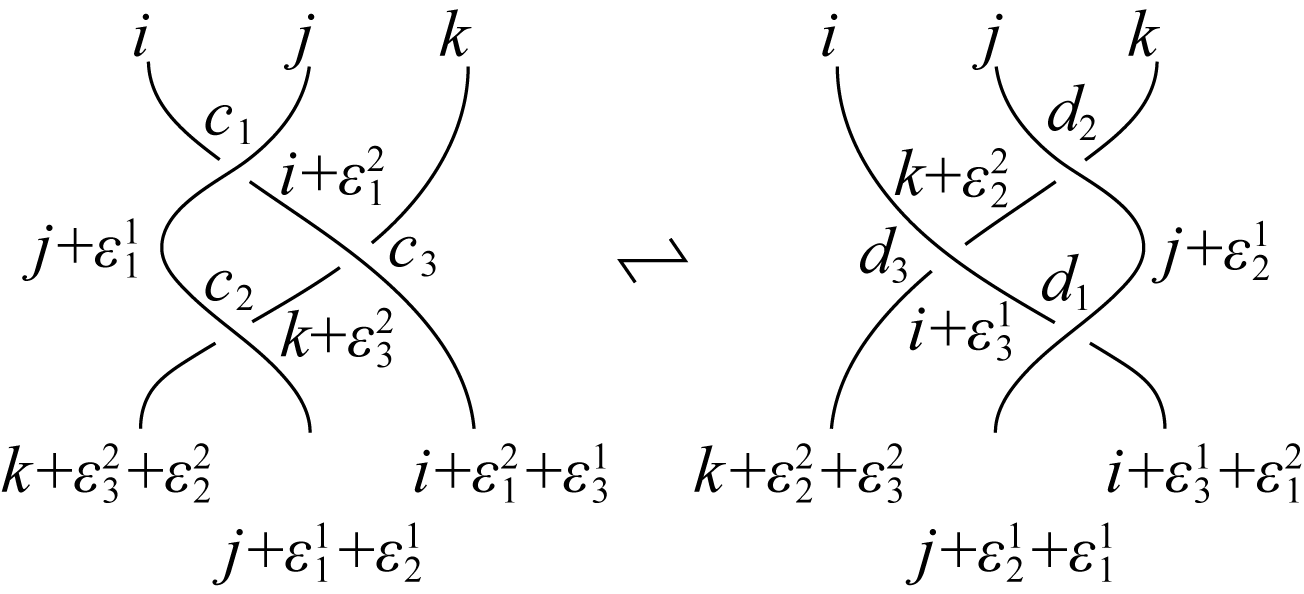}	
 \caption{Up-down labeling under R3 move}
 \label{r3}
\end{figure}

To label the V4 move, we can follow a similar procedure as used for the R3 case as shown in Figure~\ref{v4}.
\begin{figure}[h]
 \centering
 \includegraphics[width=10cm]{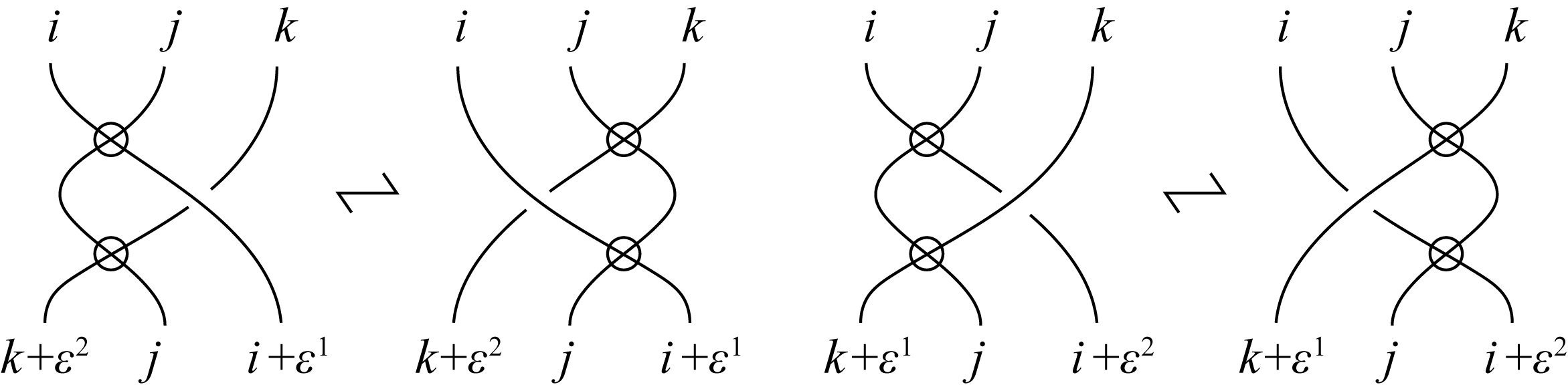}	
 \caption{Up-down labeling under V4 move}
 \label{v4}
\end{figure}

Let $c$ be the number of classical crossing in the twisted virtual braid. For the labeling under T1 and T2 move, only rule 1 of proposition~\ref{RWT} is required as shown in the Figure~\ref{t1} and Figure~\ref{t2}, respectively.
\begin{figure}[h]
	\centering
\includegraphics[width=3cm]{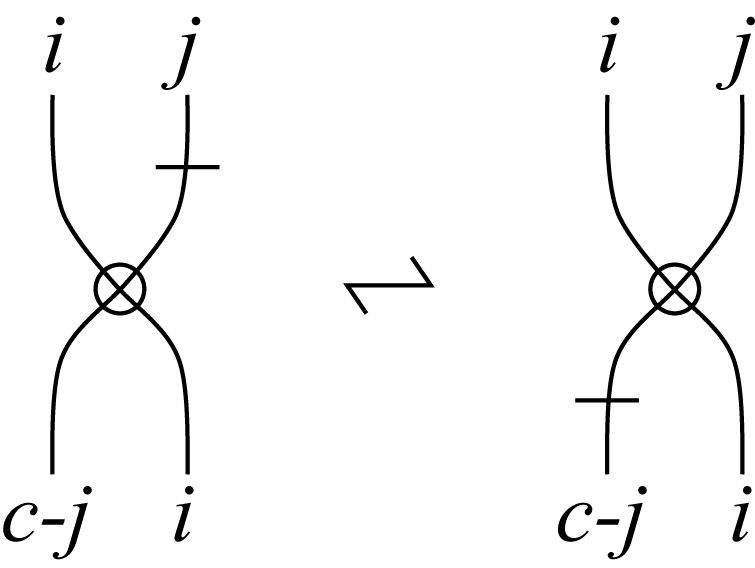}	
  \caption{Up-down labeling under T1 move}
 \label{t1}
\end{figure}

\begin{figure}[h]
	\centering
\includegraphics[width=2.5cm]{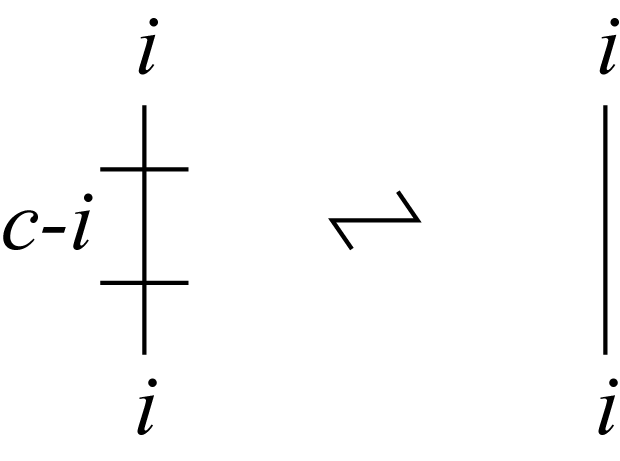}	
  \caption{Up-down labeling under T2 move}
 \label{t2}
\end{figure}
For T3 move, see Figure \ref{t3}. For the braid on the left-hand side, the first strand has an over-crossing at the crossing $d$ and the second strand has an under-crossing. Denote by $\varepsilon^1$, $\varepsilon^2$ the change at $d$ on the first, second strands, respectively. 
Note that $\varepsilon^1$ is either $0$ or $1$ and $\varepsilon^2$ is either $0$ or $-1$ in Figure~\ref{t3}. 
For the braid on the right-hand side, the first strand has an under-crossing at $d$ and the change at $d$ is $- \varepsilon^1$. The second strand has an over-crossing at $d$ and the change at $d$ is $-\varepsilon^2$. With the changes at bars, the two braids have the same labeling at the bottom. 
\begin{figure}[h]
	\centering
\includegraphics[width=5cm]{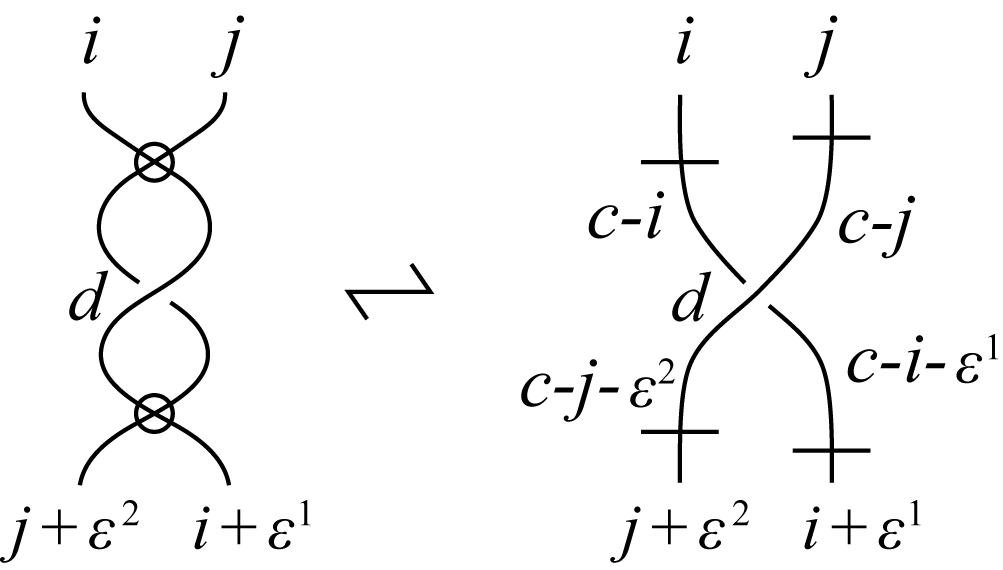}	
  \caption{Up-down labeling under T3 move}
 \label{t3}
\end{figure}
\end{proof}
\begin{corollary}\label{r2indicator12}
     If $\beta$ and $\beta'$ are equivalent virtual braids upto all virtual and classical braid moves except $R2$, then the map $f_\beta(x_1,x_2,\ldots,x_n)=f_{\beta'}(x_1,x_2,\ldots,x_n)$, for any $(x_1,x_2,\ldots,x_n) \in \mathbb{Z}^n$.
\end{corollary}
\begin{corollary}\label{r2indicator13}
    If $\beta$ and $\beta'$ are equivalent braids upto classical braid move $R3$, then the map $f_\beta(x_1,x_2,\ldots,x_n)=f_{\beta'}(x_1,x_2,\ldots,x_n)$, for any $(x_1,x_2,\ldots,x_n) \in \mathbb{Z}^n$.
\end{corollary}
\subsection{R2 indicator for twisted virtual braids}
By applying Theorem~\ref{r2indicator1}, we can establish the following theorem.

\begin{theorem}\label{r2indicator2}
    Let $\beta$ and $\beta'$ be twisted virtual braid diagrams which represent the same twisted virtual braid. If $f_\beta\neq f_{\beta'}$, then any sequence of R-moves between $\beta$ and $\beta'$ includes at least one R2 move.
\end{theorem}
\begin{corollary}\label{r2indicator22}
    Let $\beta$ and $\beta'$ be virtual braid diagrams which represent the same virtual braid. If $f_\beta\neq f_{\beta'}$, then any sequence of virtual and classical braid moves between $\beta$ and $\beta'$ includes at least one R2 move.\footnote{In \cite{KAY}, the coloring number of $\mathbb{Z}_n$-up-down coloring for virtual link diagrams was introduced as an R2 indicator and was mentioned that it was trivial for all (classical or virtual) knot diagrams. Our function $f_\beta$ can be used as an R2 indicator for braids whose closure is either a knot or link.}
\end{corollary}
\begin{corollary}\label{r2indicator23}
    Let $\beta$ and $\beta'$ be braid diagrams which represent the same braid. If $f_\beta\neq f_{\beta'}$, then any sequence of classical braid moves between $\beta$ and $\beta'$ includes at least one R2 move.
\end{corollary}
\begin{Example}
Consider the two equivalent twisted virtual braid diagrams as shown in Figure~\ref{r2d}. Since $f_\beta(0,0)=(0,0)\neq f_{\beta'}(0,0)=(-2,2)$, they are related by at least one R2 move.
    \begin{figure}[h]
  \centering
    \includegraphics[width=3.5cm]{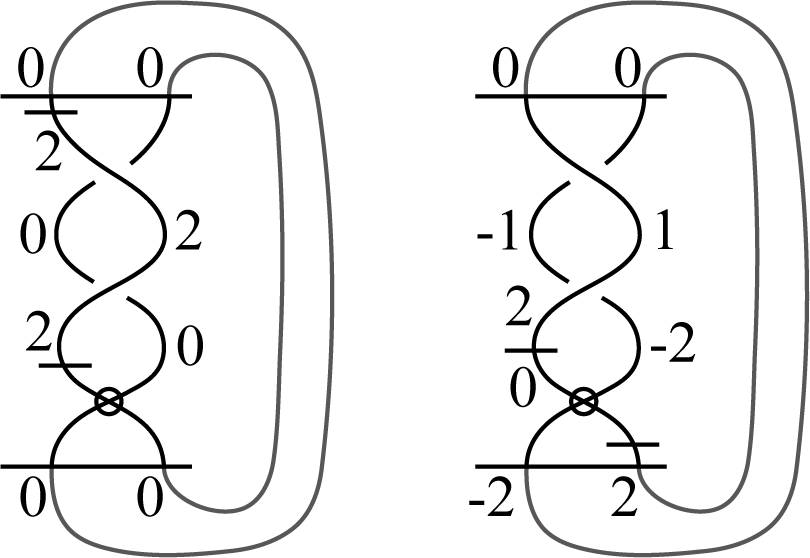}
        \caption{Two equivalent twisted virtual braid diagrams related by R2 moves}
        \label{r2d}
        \end{figure} 
\end{Example}

\begin{Example}
Consider the two equivalent virtual braid diagrams as shown in Figure~\ref{r2d2}. Since $f_\beta(0,0, 0)=(0,2,-2)\neq f_{\beta'}(0,0,0)=(-2,2,0)$, they are related by at least one R2 move.
    \begin{figure}[h]
  \centering
    \includegraphics[width=3.5cm]{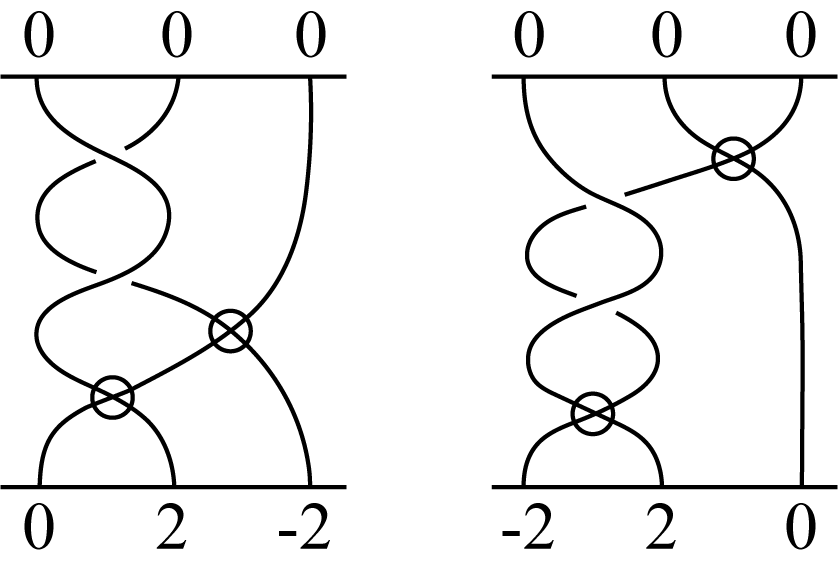}
        \caption{Two equivalent virtual braid diagrams related by R2 moves}
        \label{r2d2}
        \end{figure} 
\end{Example}
 
\section{$\mathbb{Z}_2$-labeling for twisted virtual braids}
\begin{definition}
    Let $f'_\beta$ be a map from $\mathbb{Z}^n \to \mathbb{Z}_2^n$ such that $f'_\beta=f_\beta(\text{mod }  2)$.
\end{definition}
\begin{definition}
Let $g_\beta$ be a map from $\mathbb{Z}_2^n \to \mathbb{Z}_2^n$ such that $g_\beta(x_1,x_2,\ldots,x_n)=(y_1,y_2,\ldots,y_n)$, where $x_i$ is the initial labeling of the $i^{th}$ strand from the set $\mathbb{Z}_2^n$ and $y_j$ is the labeling of $j^{th}$ strand on the bottom of $\beta$ for edges of the braids labeled using up-down labeling modulo 2 called \textit{$\mathbb{Z}_2$-labeling}. The function $g_\beta$ is denoted as \textit{$\mathbb{Z}_2$-labeling function}.
\end{definition}
\begin{Example}
In Figure \ref{weav}, four same braid diagrams with different position of a bar have the different outputs of the $\mathbb{Z}_2$-labeling functions $(0,1,1)$, $(0,0,1)$, $(1,1,0)$ and $(1,0,1)$ for the common input $(0,0,0)$. As a result, all these twisted virtual braids are mutually non-equivalent.
   \begin{figure}[h]
	\centering
\includegraphics[width=8cm]{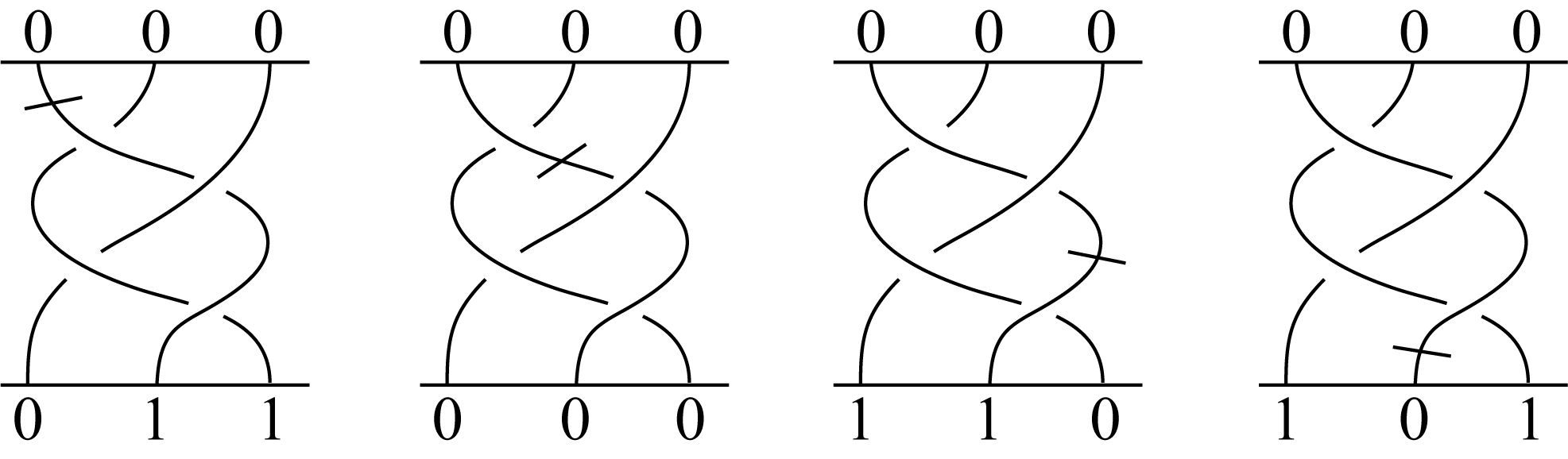}	
  \caption{Outputs of the $\mathbb{Z}_2$-labeling function}
 \label{weav}
\end{figure} 
\end{Example}

\begin{theorem}
    The functions $g_\beta$ and $f'_\beta$, will work as an invariant for twisted virtual braids.
\end{theorem}
\begin{proof}
    It is sufficient to show that it is invariant under R2 move. The Figure~\ref{r2i}, shows $g_\beta$ and $f'_\beta$ are invariant under R2 move also.
    \begin{figure}[h]
  \centering
    \includegraphics[width=4cm,height=2.5cm]{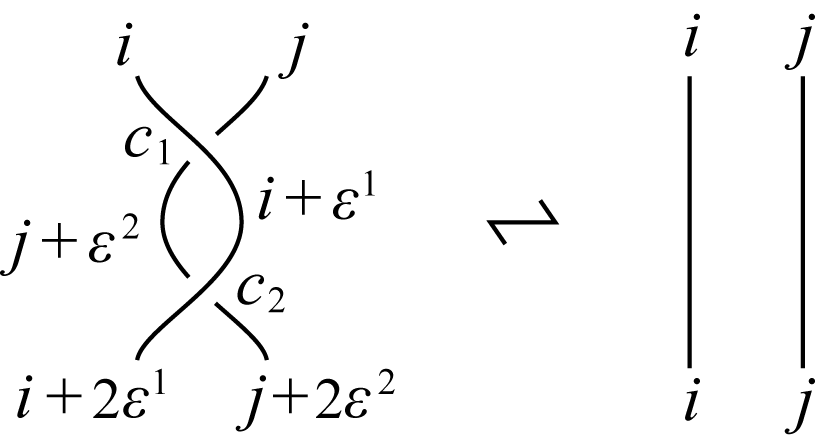}
        \caption{Labeling on R2 move}
        \label{r2i}
        \end{figure}
        
For the labeling under R2 move, see Figure \ref{r2i}. The two parts of a braid are related by an R2 move. Give the same labels, $i, j$ on the top. 
Consider crossings $c_1$ and $c_2$ in the left side of the diagram, the number of bars from $c_1$ to $c_1$ and number of bars between $c_2$ to $c_2$ are same on each strand. As a result, the labeling is conducted as depicted in the Figure~\ref{r2i}.
\end{proof}

\begin{figure}[h]
\centering
\includegraphics[width=7cm]{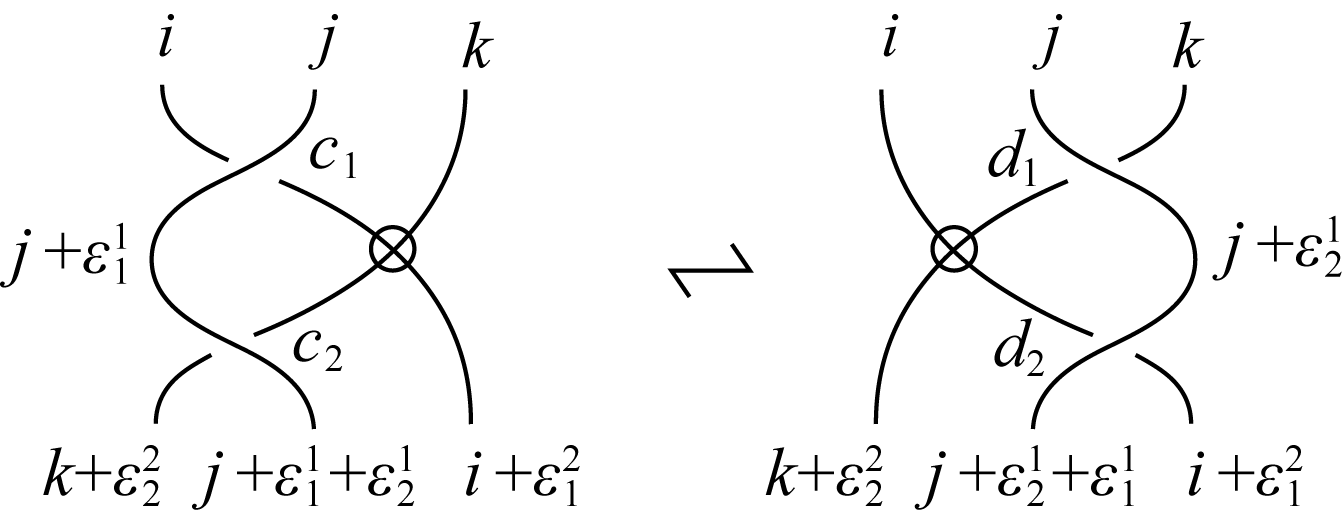}
\caption{A forbidden move on virtual knot diagram}
\label{f3}
\end{figure}

\begin{figure}[H]
\centering
\includegraphics[width=7cm]{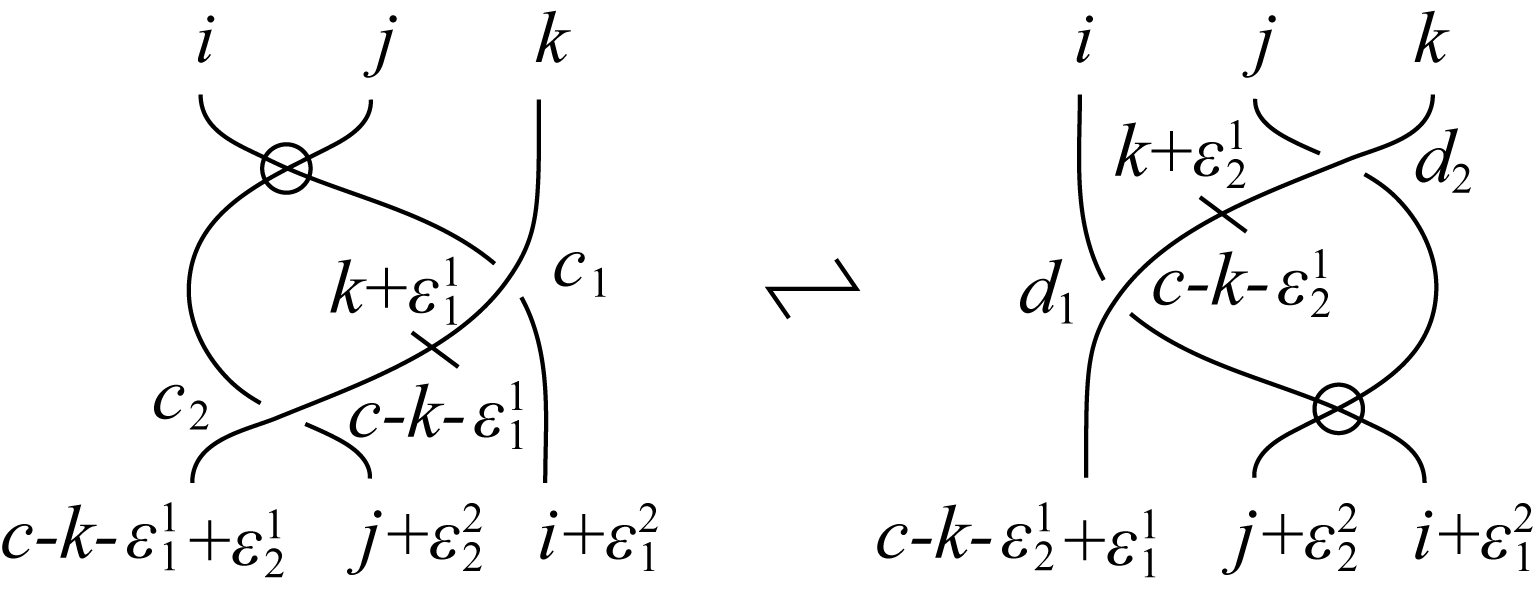}
\caption{A forbidden move on twisted knot diagram}
\label{f2}
\end{figure}

\begin{remark}
    If $\beta_0$ is a trivial twisted virtual braid with $n$ strands, then $g_{\beta_0}(x_1,x_2,\ldots,x_n)=(x_1,x_2,\ldots,x_n)$, for all $(x_1,x_2,\ldots,x_n) \in \mathbb{Z}_2^n$. Therefore, for any twisted virtual braid $\beta$, if $g_{\beta}(x_1,x_2,\ldots,x_n)\neq (x_1,x_2,\ldots,x_n)$, for some $(x_1,x_2,\ldots,x_n) \in \mathbb{Z}_2^n$, then $\beta$ is a non-trivial twisted virtual braid.
\end{remark}
Using the function $g_\beta$, we can define a one variable polynomial for $\beta$ as follows: 
\begin{definition}
    Let $g_\beta(x_1,x_2,\ldots,x_n)=(y_1,y_2,\ldots,y_n)$ satisfy the following relation: $y_{p_{\beta}(i)}=L_{x_i}(x_i)$, where $p_\beta$ be the permutation of $\beta$. Define a one variable polynomial corresponding to $\beta$, $F_\beta(x)$ in $\mathbb{Z}_2[x]$ as given below:
    \[F_\beta(x)=\prod_{i=1}^n L_{x_i}(x).\] It is an invariant for twisted virtual braids, named as $\mathbb{Z}_2$-polynomial.
\end{definition}
\begin{Example}
    In Figure \ref{weav}, four twisted virtual braid diagrams are given. The $\mathbb{Z}_2$-polynomial of second twisted virtual braid from the left is $x^3$ and for others, it is $x(x+1)^2$. This implies that $g_\beta$ is a stronger invariant than $\mathbb{Z}_2$-polynomial.
\end{Example}
In twisted knot theory, there are three independent forbidden moves, as described in reference~\cite{SQ}. The up-down labeling remains invariant under the forbidden move allowed in welded knot theory depicted in Figure \ref{f3} but it is not invariant under the second forbidden move depicted in Figure~\ref{f2}.


\begin{Example}
In the Figure~\ref{exf}, the two diagrams are related by forbidden move and extended Reidemeister moves and the labeling output is different for both, i.e. $g_\beta(0,0,0)=(1,0,1)\neq g_{\beta'}(0,0,0)=(1,1,0)$. In Figure~\ref{exf2}, a sequence of R-moves and a forbidden move is depicted with output of $(0,0,0)$ under $\mathbb{Z}_2$-labeling function.
\begin{figure}[h]
\centering
\includegraphics[width=5cm]{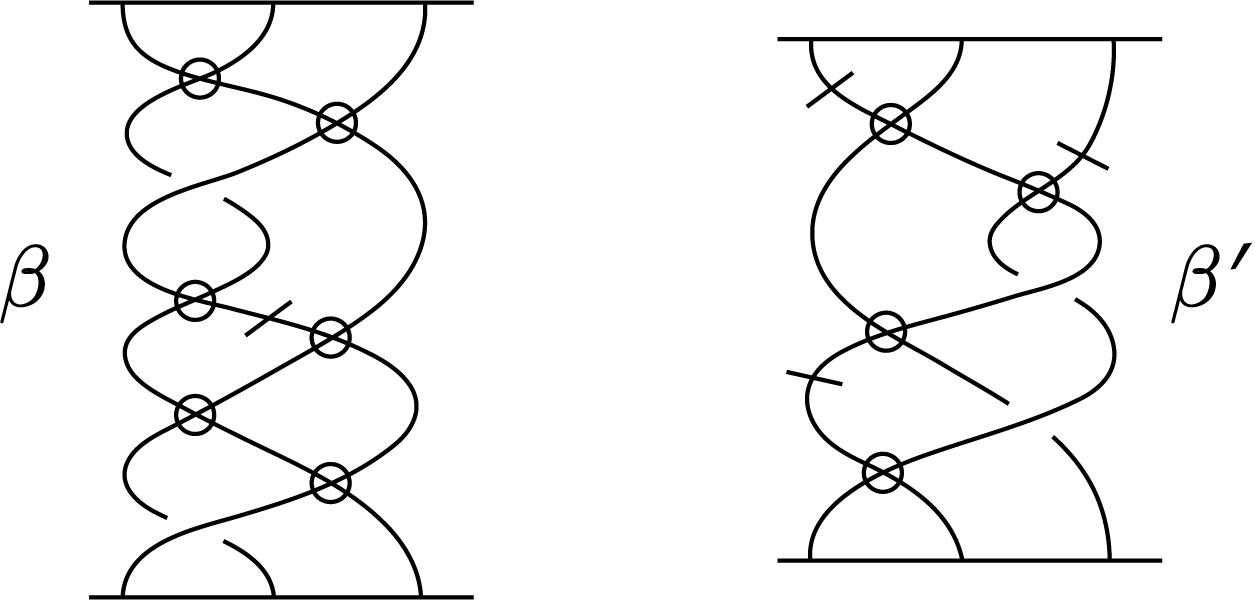}
\caption{A pair of twisted virtual braid diagrams which are transformed into each other by extended Reidemeister moves and a forbidden move}
\label{exf}
\end{figure} 

\begin{figure}[H]
\centering
\includegraphics[width=12cm]{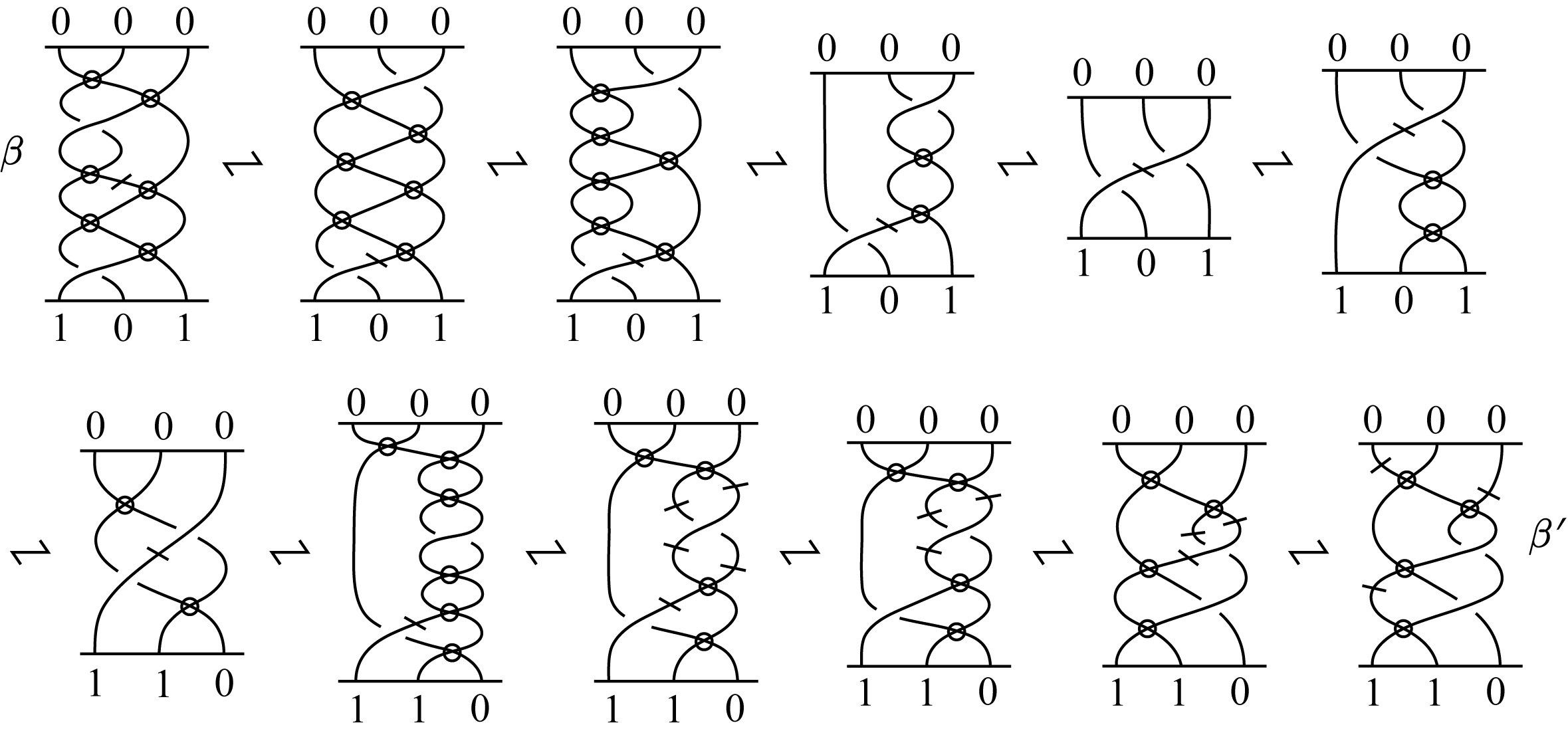}
\caption{A sequence of moves between $\beta$ and $\beta'$, including a forbidden move}
\label{exf2}
\end{figure} 
\end{Example}

\section*{Conclusion}
In this paper, our primary focus is to explore the labeling invariant for twisted virtual braids. 
We also plan to investigate the application of this labeling concept to twisted link diagrams in future studies, which represents an intriguing direction for research.
Furthermore, a compelling challenge lies in developing invariants based on the warping degree of twisted knot diagrams. 
This challenge motivates further exploration into understanding and characterizing twisted knot structures using quantitative measures like the warping degree.

\section*{Acknowledgements}
The first author would like to thank the University Grants Commission(UGC), India, for Research Fellowship with NTA Ref.No.191620008047. 
The second author's work was partially supported by JSPS KAKENHI Grant Number JP21K03263. The third author acknowledges the support given by SERB research project(MATRICS) with F.No.MTR/2021/000394 and by the NBHM, Government of India under grant-in-aid with F.No.02011/2/20223 NBHM(R.P.)/ R\&D II/970.

\noindent Komal Negi \\
Department of Mathematics, Indian Institute of Technology Ropar, Punjab, India.
\begin{verbatim} komal.20maz0004@iitrpr.ac.in \end{verbatim} 

\noindent Ayaka Shimizu \\ 
Osaka Central Advanced Mathematical Institute (OCAMI), Osaka Metropolitan University, 
3-3-138, Sugimoto, Osaka, 558-8585, Japan. 
\begin{verbatim} shimizu1984@gmail.com \end{verbatim}

\noindent Madeti Prabhakar \\
Department of Mathematics, Indian Institute of Technology Ropar, Punjab, India.
\begin{verbatim} prabhakar@iitrpr.ac.in \end{verbatim} 

\end{document}